\documentstyle[12pt]{article}

\def\C{{\rm C  \! \! \! l \ }}
\def\R{{\rm I   \! R \ }}

\begin{document}
\begin{flushright} 
{ q-alg/0102071}
\end{flushright}
\title{POSSIBLE CONTRACTIONS OF QUANTUM ORTHOGONAL GROUPS}
\author{N.A.Gromov, I.V. Kostyakov, V.V. Kuratov \\
Department of Mathematics, \\
Syktyvkar Branch of IMM UrD RAS,\\
e-mail:gromov@dm.komisc.ru} 

\maketitle

\begin{abstract}

Possible contractions of quantum orthogonal groups which
correspond to different
choices of primitive elements of
Hopf algebra  are considered
 and all allowed contractions in Cayley--Klein
scheme are obtained. Quantum deformations of  kinematical
groups have been investigated and have shown that  quantum
analog of  (complex) Galilei group $G(1,3)$ do not exist in our scheme.
\end{abstract}

\section{Introduction}

Contraction of Lie groups (algebras) is a method of obtaining new Lie
groups (algebras) from some initial ones with the help of passage to the
limit \cite{IW}. One may define contraction of algebraic structure $(M,*)$
as the map $\phi_{\epsilon} : (M,*) \rightarrow (N,*')$,
where $(N,*)$ is  algebraic structure of the same type,
isomorphic to $(M,*)$ for  $ \epsilon \neq 0 $ and nonisomorphic to
the initial one for $ \epsilon = 0. $
Except for Lie group (algebra) contractions, graded contractions
\cite{P,MP} are known, which preserve the grading of Lie algebra.
Under contractions of bialgebra \cite{VG}
Lie algebra structure and cocommutator are conserved.
Hopf algebra (or quantum group) contractions are introduced in
such a way \cite{C1,C2}, that in the limit $ \epsilon \rightarrow 0 $
a new expressions for coproduct, counit and antipode
are consistent with  Hopf algebra axioms.

Contractions as a passage to  limit are corresponded with a
physical intuition. At the same time it is desirable
to investigate contractions of an algebraic structures
with the help of pure algebraic tools. It is possible for
classical and quantum  groups and algebras if one take into
consideration Pimenov algebra $ {\bf D}(\iota) $ with nilpotent
commutative generators $ \cite{2}$.

In present paper contractions of  quantum orthogonal
groups are studed and the groups under consideration are regarded
according to \cite{Fad-89} as an algebra of noncommutative functions
but with nilpotent generators. From the contraction viewpoint
Hopf algebra structure of quantum orthogonal group is more rigid
as compared with the group one. Possible contractions are essentially
depend on the choice of primitive elements of  Hopf algebra.
We have regarded all variants of such choise for quantum orthogonal
group $ SO_q(N)$ and for each variant have find all admissible contractions
in Cayley-Klein scheme.

\section{Orthogonal Cayley-Klein groups}

Let us define {\it Pimenov algebra} $ {\bf D}_n(\iota; {\C}) $
as an associative algebra with unit over complex number field
and with nilpotent commutative generators
$ {\iota}_k, \ {\iota}_k^2=0,$
$ {\iota}_k{\iota}_m={\iota}_m{\iota}_k \not =0, \  k \neq m,
\  k,m=1, \ldots, n. $
The general element of ${\bf D}_n(\iota; {\C}) $
is in the form
\begin{equation}
d=d_0+\sum^{n}_{p=1}\sum_{k_1< \ldots < k_p}d_{k_1\ldots k_p}
{\iota}_{k_1} \ldots {\iota}_{k_p}, \quad  d_0,d_{k_1 \ldots k_p} \in {\C}.
\label{3.1}
\end{equation}
It is possible to  define  the division of nilpotent
generator  $ {\iota}_k $ by itself, namely:
$ {\iota}_k/{\iota}_k=1, \ k=1, \ldots, n. $
Let us stress that the division of different nilpotent generators
$ {\iota}_k/{\iota}_p, \ k \not = p, $ as well as the division of
complex number by nilpotent generators
$ a/{\iota}_k, \ a \in {\C} $  are not define.

Let $ SO(N;{\C}) $ be an orthogonal matrix group.
Its elements are matrices $ A=(a_{kp})\in {M_N (\C)},$
$ A^t=A^{-1} $ and under the action
$ y'=Ay $ on  vectors $ y $ of  complex vector space $O_N$
the quadratic form $ y^{t}y=\sum\nolimits_{k=1}^{N}y^2_k$ is preserved,
where $ y_k$ are Cartesian components of $y$.
Sometimes it is convenient to regard an orthogonal group in so-called
"symplectic" basis. Transformation from Cartesian to symplectic basis
$ x=Dy $ is made by matrix $D$, which is a solution of equation
     \begin{equation}
 D^{t}C_{0}D=I,
     \label{1-1}
     \end{equation}
where $C_{0}\in M_{N},\; (C_{0})_{ik}= \delta_{i,k'},\; k'=N+1-k$
Equation (\ref{1-1}) has many solutions, take one of them, namely
\begin{equation}
D=\frac{1}{\sqrt{2}}
\left ( \begin{array}{ccc}
      I & 0 &  -i{\tilde C_0} \\
      0 & \sqrt{2} &  0 \\
      {\tilde C_0} & 0  &  iI
      \end{array} \right ),    \       N=2n+1,
\label{1-3}
\end{equation}
where $n \times n $ matrix $ {\tilde C_0} $ is similar to
$C_0.$ For $N=2n$ the matrix $D$ is given by (\ref{1-3}) but without
middle  column and row.
Matrices $ B $ of $ SO(N;{\C}) $ in symplectic basis are obtained from
$ A $ by similarity transformation $ B=DAD^{-1} $ and are
subject of orthogonality relations $ B^{t}C_{0}B=C_{0}. $
The quadratic form $ x^{t}C_{0}x $
is invariant under the action $x'=Bx$.

{\it Complex orthogonal Cayley-Klein group} $ SO(N;j;\C) $
is defined as the group of transformations
$ {\xi}'(j)=A(j)\xi(j)$ of complex vector space $O_N(j)$
with Cartesian coordinates
$ \xi^t (j)= (\xi_1,  (1,2)  \xi_2,  \ldots,(1,N)\xi_N)^t,
 $
which preserve the quadratic form
$
inv(j) =\xi^t(j)\xi(j) = \xi^2_1 + \sum^N_{k=2}({1,k})^2{\xi}^2_k,
$
where $ j=(j_1, \ldots, j_{N-1}) $, each parameter $ j_k $
takes {\it two} values:
$  j_r=1,{\iota}_r, \ r=1, \ldots, N-1, $
and
\begin{equation}
(\mu,\nu)=\prod^{max(\mu,\nu)-1}_{l=min(\mu,\nu)} j_l, \quad (\mu,\mu)=1.
\label{3.3}
\end{equation}
Let us stress, that Cartesian coordinates
of $  O_N(j) $ are special elements of Pimenov algebra
${\bf D}_{N-1}(j; \C). $
 Cayley-Klein group $ SO(N;j;{\C}) $ in turn may be
realised as an matrix group whose elements are taken from algebra
$ {\bf D}_{N-1}(j;{\C}) $
and consist of the $ N \times N $ matrices $A(j)$
$
(A(j))_{kp}=(k,p)a_{kp}, \ a_{kp} \in \C.
$
Matrices $A(j)$ are subject of the additional $j$-orthogonality relations
$
A(j)A^t(j) = A^t(j)A(j)=I.
$

The passage to the symplectic description is made by  matrices,
which are solutions of equations (\ref{1-1}).
Let us regard the matrix $ D_{\sigma}=DV_{\sigma},$ where
$V_{\sigma} \in M_{N},$ $(V_{\sigma})_{ik}= \delta_{\sigma_{i},k},$
and $ \sigma \in S(N) $ is a permutation of the $ N$-th order.
It is easy to verify that $ D_{\sigma}$
is again a solution of equation (\ref{1-1}).
Then in symplectic basis the orthogonal Cayley Klein group
$SO(N;j;{\C})$ is described by the matrices
$
 B_{\sigma}(j)=D_{\sigma}A(j)D^{-1}_{\sigma}
 $
with the additional relations of $ j$-orthogonality
$
   B_{\sigma}(j)C_{0}B^{t}_{\sigma}(j)=
B^{t}_{\sigma}(j)C_{0}B_{\sigma}(j)=C_{0}.
$

It should be noted that for  orthogonal groups $(j=1)$
the use of different matrices $D_{\sigma}$ has no sense because of all
Cartesian coordinates of $O_N$ are equivalent
up to a choice of its enumerations.
The different situation is for Cayley-Klein groups $(j \not =1).$
Cartesian coordinates
$(1,k)\xi_{k}, \ k=1,\dots,N$ for nilpotent values of some or all
parameters $j_k$ are different elements of the algebra $D_{N-1}(j;\C),$
therefore the same group $SO(N;j;\C)$ may be realized by  matrices
$ B_{\sigma} $  with a {\it different} disposition of nilpotent
generators among their elements.
Matrix elements of $ B_{\sigma}(j) $ are as follows
\begin{equation}
     \begin{array}{l}
(B_{\sigma})_{n+1,n+1}= b_{n+1,n+1},  \\
(B_{\sigma})_{kk}=b_{kk}+i\tilde{b}_{kk}(\sigma_{k},
\sigma_{k'}), \quad
(B_{\sigma})_{k'k'}=b_{kk}-i\tilde{b}_{kk}(\sigma_{k},
\sigma_{k'}), \\
(B_{\sigma})_{kk'}=b_{k'k}-i\tilde{b}_{k'k}(\sigma_{k},
\sigma_{k'}), \quad
(B_{\sigma})_{k'k}=b_{k'k}+i\tilde{b}_{k'k}(\sigma_{k},
\sigma_{k'}), \\
(B_{\sigma})_{k,n+1}=b_{k,n+1}(\sigma_k, \sigma_{n+1}) -
i\tilde{b}_{k,n+1}(\sigma_{n+1}, \sigma_{k'}), \\
(B_{\sigma})_{k',n+1}=b_{k,n+1}(\sigma_k, \sigma_{n+1}) +
i\tilde{b}_{k,n+1}(\sigma_{n+1}, \sigma_{k'}), \\
(B_{\sigma})_{n+1,k}=b_{n+1,k}(\sigma_k, \sigma_{n+1}) +
i\tilde{b}_{n+1,k}(\sigma_{n+1}, \sigma_{k'}), \\
(B_{\sigma})_{n+1,k'}=b_{n+1,k}(\sigma_k, \sigma_{n+1}) -
i\tilde{b}_{n+1,k}(\sigma_{n+1}, \sigma_{k'}), \; k \neq p, \\
(B_{\sigma})_{kp}=b_{kp}(\sigma_k, \sigma_{p})+
b_{kp}'(\sigma_{k'},\sigma_{p'})+
i\tilde{b}_{kp}(\sigma_{k}, \sigma_{p'})-
i\tilde{b}_{kp}'(\sigma_{k'}, \sigma_{p}),  \\
(B_{\sigma})_{kp'}=b_{kp}(\sigma_k, \sigma_{p})-
b_{kp}'(\sigma_{k'},\sigma_{p'})-
i\tilde{b}_{kp}(\sigma_{k}, \sigma_{p'})-
i\tilde{b}_{kp}'(\sigma_{k'}, \sigma_{p}),  \\
(B_{\sigma})_{k'p}=b_{kp}(\sigma_k, \sigma_{p})-
b_{kp}'(\sigma_{k'},\sigma_{p'})+
i\tilde{b}_{kp}(\sigma_{k}, \sigma_{p'})+
i\tilde{b}_{kp}'(\sigma_{k'}, \sigma_{p}),  \\
(B_{\sigma})_{k'p'}=b_{kp}(\sigma_k, \sigma_{p})+
b_{kp}'(\sigma_{k'},\sigma_{p'})-
i\tilde{b}_{kp}(\sigma_{k}, \sigma_{p'})+
i\tilde{b}_{kp}'(\sigma_{k'}, \sigma_{p}) . \\
\end{array}
   \label{2-1}
\end{equation}
Here $ b,  b', \tilde{b}, \tilde{b}' \in \C$
may be easily expressed by  matrix elements of $A. $

\section{ Contractions of quantum orthogonal groups. }

\subsection{Formal definition of the quantum group $ SO_v(N;j;\sigma) $ }

 The starting point of the definition of quantum groups \cite{Fad-89}
 is an algebra $ \C \langle T_{ik} \rangle $
of noncommutative polynomials of $ N^2 $ variables.
We start with an algebra $ {\bf D} \langle (T_{\sigma})_{ik} \rangle $
of noncommutative polynomials of  $ N^2 $ variables, which are an elements of
the direct product $ {\bf D}_{N-1}(j)\otimes {\C}\langle t_{ik} \rangle. $
More precisely the elements $  (T_{\sigma})_{ik} $ are obtained
from the elements $ \left( B_{\sigma}(j) \right)_{ik}$ of 
equations~(\ref{2-1})
by the replacement of commutative variables $ b,b',\tilde{b},\tilde{b}' $
with the noncommutative variables $ t,t',\tau,\tau',$ respectively.
One introduce additionally the transformation of the deformation parameters
$q=e^z $ as follows:
$ z=Jv, $
where $ v $ is a new deformation parameter and  $ J $ is some
product of  parameters $ j $ for the present unknown.
Let $ \tilde R_v(j), C(j) $ are matrices which are obtained from
corresponding matrices of \cite{Fad-89}
by the replecement of deformation parameter $ z $ with $ Jv: $
\begin{equation}
R_{v}(j)=R_{q}(z \rightarrow Jv), \quad
C(j)=C(z \rightarrow Jv).
\label{2-4}
\end{equation}
The commutation relations of the generators
$ T_{\sigma}(j) $ are defined by
\begin{equation}
 R_v(j)T_1(j)T_2(j)=T_2(j)T_1(j)R_v(j),
\label{21}
\end{equation}
where $ T_1(j)=T_{\sigma}(j) \otimes I, \
T_2(j)=I \otimes T_{\sigma}(j)   $
and the additional relations of $ (v,j)$-orthogonality
     \begin{equation}
   T_{\sigma}(j)C(j)T^{t}_{\sigma}(j) =
T^{t}_{\sigma}(j)C(j)T_{\sigma}(j) = C(j)
    \label{22}
    \end{equation}
are imposed.

{\it One define the quantum orthogonal Cayley-Klein group}
$ SO_{v}(N;j;\sigma) $ {\it as the quotient algebra of}
 $ {\bf D} \langle (T_{\sigma})_{ik} \rangle $
{\it by relations} (\ref{21}), (\ref{22}).
Formally $ SO_{v}(N;j;\sigma)$ is a Hopf algebra with the following coproduct
$ \Delta, $ counit $ \epsilon $ and antipode $ S: $
\begin{equation}
\Delta T_{\sigma}(j)=T_{\sigma}(j) \dot {\otimes}T_{\sigma}(j),
\quad \epsilon (T_{\sigma}(j))=I, \quad
S(T_{\sigma}(j))=C(j)T^{t}_{\sigma}(j)C^{-1}(j).
\label{25}
\end{equation}
As far as only second diagonal elements of the matrix C are
different from zero and for $ q=1 $ this matrix is equal to
$ C_{0}, $ then we have the symplectic description of
$ SO_{q}(N).$

\subsection{Allowed contractions of $ SO_v(N;j;\sigma) $ }

The formal definition of
$ SO_{v}(N;j;\sigma) $ should be a real definition of quantum group,
if the proposed construction is a consistent Hopf algebra structure
under nilpotent values of some or all parameters $j$.
Counit $ \epsilon(t_{n+1,n+1})=1, \; \epsilon(t_{kk})=1, \;
k=1,\ldots,n $ and  $\epsilon(t)=\epsilon(\tau)=0 $
for the rest generators
do not restrict the values of  $ j. $
Parameters $j$ are arranged in the expressions for coproduct
$ \Delta $  exactly as in matrix product of
$ B_{\sigma}(j), $ and as far as the last ones form the group
$ SO(N;j;\C)$ for any values of $j,$ then no restrictions
follow from the coproduct. Different situation is with the
antipode $ S. $  Really, for elements
$
(T_{\sigma})_{k'k}=t_{k'k} + i\tau_{k'k}(\sigma_k,\sigma_{k'}),
  \quad  k=1,\ldots,n,
$
 antipode is obtained as
\begin{equation}
S((T_{\sigma})_{k'k})=(T_{\sigma})_{k'k}\cdot e^{2J\rho_kv}
\label{15-2}
\end{equation}
and depend both on $ \rho_k$ and for the present undetermined
factor $J.$ Antipode is an antihomomorphism of Hopf algebra and
therefore have to transform  $ T_{\sigma}(j) $ to a matrix
with the same distribution of the nilpotent parameters $j$ in its elements,
i.e. the right and the left parts of equation (\ref{15-2})
must be identical elements of
$ {\bf D}_{N-1}(j)\otimes {\C}\langle t_{ik} \rangle. $
For $J=1$ this condition is holds for any values of the parameters $j.$
The case $ J \neq 1$ require additional discussion.

Next condition which must be taken into account is the
$(v,j)$-orthogonality relations (\ref{22}).
We require that the number of equations in (\ref{22})
 is not changed as compared
with the initial quantum group. It is possible when nilpotent generators
are appeared in equation (\ref{22}) either with the powers greater or
equal two (and then the corresponding terms are equal to zero) or
as homogeneous multipliers.
Taking into account all these arguments and using the explicit
expressions for antipode and $(v,j)$-orthogonality we can find
possible contractions of quantum orthogonal groups, which
are described by the following theorems.

{ \bf Theorem 1.}   { \it
If the deformation parameter is not transformed $J=1$,
then the following maximal $n$-dimensional contraction of the
orthogonal quantum group  $ SO_v(N;j;\sigma), \ N=2n+1 $ is allowed:
$
 j_{2s}=\iota_{2s}, \; s=1,\ldots,m, \;
j_{2r+1}=\iota_{2r+1},  \; r=m,\ldots,n-1, \; 0 \leq m \leq n,
  $
for example, for permutation $\sigma $:
$
\sigma_{n+1}=2m+1, \;
\sigma_{s}=2s-1, \; \sigma_{s'}=2s, \; s=1,\ldots,m, \;
\sigma_{r}=2r, \; \sigma_{r'}=2r+1, \; r=m+1,\ldots,n.
$ }

{ \bf Theorem 2.}   { \it
If the deformation parameter is not transformed $J=1$,
then the following maximal n-dimensional contraction of the quantum orthogonal
group  $ SO_v(N;j;\sigma), \ N=2n $ is allowed:
$
 j_{2s}=\iota_{2s}, \; s=1,\ldots,m-1, \;
j_{2p-1}=\iota_{2p-1},  \; p=m,\ldots,u, \;
 j_{2r}=\iota_{2r},  \; r=u,\ldots,n-1, \;
1 \leq m \leq u \leq n,
 $
for example, for permutation $ \sigma $:
$ \sigma_{n}=2m-1, \;\sigma_{n'}=2u, \;
\sigma_{s}=2s-1, \; \sigma_{s'}=2s, \; s=1,\ldots,m-1, \;
\sigma_{p}=2p, \; \sigma_{p'}=2p+1, \; p=m,\ldots,u-1, \;
\sigma_{r}=2r+1, \; \sigma_{r'}=2r, \; r=u,\ldots,n-1.
$ }

{\bf Remark 1.}
It should be noted that as $ \sigma $ may be taken any permutation
with the properties $ (\sigma_{k},\sigma_{k'})=1, \; k=1,\ldots,n $
(or $ n-1$).

{\bf Remark 2.}
Admissible contractions for number of parameters  $ j_k $
less then $n$ are obtained from theorems 1 and 2
by setting part of $ j_{2s},  j_{2p-1},  j_{2r},$ $ j_{2r+1} $ equal to one.

We return to the antipode  (\ref{15-2}) for $ J \neq 1.$
As far as $ \rho_{n+1}=0 $ for $ N=2n+1,$ and $ \rho_{n} = \rho_{n'}=0 $
for $N=2n,$ we shall regard these two cases separately.

{ \bf Theorem 3.}   { \it
If the deformation parameter is  transformed $ (J \not =1),$
then the following contractions of the quantum orthogonal
group  $ SO_v(N;j;\sigma),$ $ N=2n+1 $ are allowed:

1. For  $ J=j_{n+1}, $
a) $ j_{n+1}=\iota_{n+1},$  if
$ 1<\sigma_{n+1}<n+1; \quad $
b) $ j_{n+1}=\iota_{n+1}, \; j_{1}=1,\iota_{1},$
 if $ \sigma_{n+1} = 1.$

2. For  $ J=j_{n}, $
a) $ j_{n}=\iota_{n}, $  if
 $n+1<\sigma_{n+1}<2n+1; \quad $
b) $ j_{n}=\iota_{n}, \;  j_{2n}=1,\iota_{2n},$
 if $ \sigma_{n+1} =2n+1. $

3. For  $  J=j_nj_{n+1}, \;$
$j_{n}=1,\iota_{n},\; j_{n+1}=1,\iota_{n+1}, $
 if  $ \sigma_{n+1}=n+1. $
}

{ \bf Theorem 4.}   { \it
     If the deformation parameter is transformed $(J\neq1),$
then the following contractions of the quantum orthogonal group
$SO_v(N;j;\sigma), $ $ N=2n $ are allowed:

1. For $ J=j_{n}, $
a) $ j_{n}=\iota_n, $  if
$ \sigma_{n}>1,\;  \sigma_{n'}<2n; $
b) $ j_{n}=\iota_n, \; j_{1}=1,\iota_{1}, $  if
$ \sigma_{n}=1,\;  \sigma_{n'}< 2n; $
c) $ j_{n}=\iota_n, \; j_{2n-1}=1,\iota_{2n-1}, $  if
$\sigma_n > 1,\; \sigma_{n'}= 2n; $
d) $ j_{n}=\iota_n, \; j_{1}=1,\iota_{1},
\; j_{2n-1}=1,\iota_{2n-1}, $
if $\sigma_n=1,\;  \sigma_{n'}= 2n.  $

2. For $ J=j_{n-1}. $
a) $ j_{n-1}=\iota_{n-1}, $  if
$ \sigma_{n'}< 2n; $
b) $ j_{n-1}=\iota_{n-1}, \; j_{2n-1}=1,\iota_{2n-1}, $
 if $ \sigma_n<2n-1, \; \sigma_{n'}=2n; $
c) $ j_{n-1}=\iota_{n-1}, \; j_{2n-2}=1,\iota_{2n-2}, \;
 j_{2n-1}=1,\iota_{2n-1},  $
 if $\sigma_n=2n-1,\;  \sigma_{n'}= 2n.  $

3. For  $ J=j_{n+1}. $
a) $ j_{n+1}=\iota_{n+1}, $
 if $ \sigma_{n}>1; $
b) $ j_{n+1}=\iota_{n+1}, \; j_{1}=1,\iota_{1}, $
 if $ \sigma_{n}=1, \; \sigma_{n'}>2; $
c) $ j_{n+1}=\iota_{n+1}, \; j_{1}=1,\iota_{1}, \;
j_{2}=1,\iota_{2},  $  if $ \sigma_{n}=1,\;  \sigma_{n'}= 2.  $

4. For  $ J=j_{n-1}j_{n}. $
a) $ j_{n-1}=\iota_{n-1}, \; j_{n}=\iota_{n}, $
 if $ \sigma_{n'}<2n; $
b) $ j_{n-1}=\iota_{n-1}, \; j_{n}=\iota_n, \; j_{2n-1}=1,\iota_{2n-1}, $
 if $ \sigma_{n'}=2n.  $

5. For  $ J=j_{n}j_{n+1}. $
a) $ j_{n}=\iota_{n}, \; j_{n+1}=\iota_{n+1}, $
 if $\sigma_{n}>1; $
b) $ j_{1}=1,\iota_{1}, \; j_{n}=\iota_{n}, \; j_{n+1}=\iota_{n+1}, $
 if $ \sigma_{n}=1. $

6. For  $ J=j_{n-1}j_{n}j_{n+1}. $
a) $ j_{n-1}=1,\iota_{n-1}, \; j_{n}=1,\iota_{n}, \; j_{n+1}=1,\iota_{n+1}, $
 if
 $ \sigma_{n}=n,\; \sigma_{n'}=n+1. $
    }

Hopf algebra  $ SO_q(N;j;\sigma),\ N=2n+1$ has $n$
primitive elements which correspond to $n$
 diagonal $ 2\times2 $ submatricies:
diag$((B_{\sigma})_{kk},(B_{\sigma})_{k'k'}) =
$ diag$(b_{kk}+i\tilde{b}_{kk}(\sigma_{k},\sigma_{k'}),
     b_{kk}-i\tilde{b}_{kk}(\sigma_{k},\sigma_{k'})),
     k=1,\ldots,n, $
see (\ref{2-1}).
If the deformation parameter $z$ is fixed
$(J=1)$ under contractions, then {\it all} primitive elements
of the contracted quantum orthogonal group are corresponded to
Euclidean rotation  $ SO(2).$
If the deformation parameter is transformed $ z=\iota v, $
 then { \it all } primitive elements are correspond
to Galilei transformation  $ SO(2;j=\iota)=G(1,1).$ The same is true for the
contracted quantum groups $ SO_q(N;j;\sigma),\ N=2n.$
Let us note that contractions of quantum orthogonal algebras
with a different sets of primitive elements have been discussed in
\cite{VG}, \cite{T}.

Quantum orthogonal groups have contractions with the
same nilpotent parameters $j$ both with fixed deformation
parameter and with transformed one. For example, quantum group
$ SO_q(2n+1;j;\sigma)$ for even $ n=2p $  at
$ \sigma_{n+1}=1 $ according with theorem 1 has contraction
$ j_{n}=\iota_{n}, \; j_{n+1}=\iota_{n+1}, \; J=1$
and according with item 3 of theorem 3 has
the same two dimensional contraction,
but $  J=\iota_{n}\iota_{n+1}.$
Let us stress that the cases $ J=1$ and $J \sim  \iota $
are realized for {\it different} sets of primitive elements in
Hopf algebra.

Let permutation $\sigma$ is identical, i.e.
$\sigma_k=k,$ $\sigma_{k'}=k',$ $\sigma_{n+1}=n+1.$
It follows from theorems  1 and 2 that there are no contractions of
 $SO_q(N;j)$ with fixed deformation parameter $(J=1)$.
For $N=2n+1$ from theorem 3 we obtain three possible contractions:
$j_n=1,\iota_n,$ $j_{n+1}=1,\iota_{n+1}$
(both parameters $j_n$ and $j_{n+1}$ independently take nilpotent values)
and deformation parameters is transformed  with $J=j_nj_{n+1}.$
For $N=2n$ from theorem 4 we obtain seven admissible contractions:
$j_{n-1}=1, \iota_{n-1},$  $j_{n}=1, \iota_{n},$
$j_{n+1}=1, \iota_{n+1},$ where deformation parameter is multiplied by
$J=j_{n-1}j_nj_{n+1}.$ It should be considered in \cite{Sb-97} just
these allowed contractions.

From the contraction viewpoint  Hopf algebra structure of quantum
orthogonal group is more rigid as compared with a group one.
Cayley-Klein groups are obtained from $SO(N;j)$ for all nilpotent values
of  parameters $j_k,\ k=1,\ldots,N-1,$ whereas their quantum
deformations exist only for  some of them $(\leq[{N \over 2}])$.
It should be note that among the contracted for equal number of
parameters $j$ quantum orthogonal groups may be isomorphic as Hopf
algebra quantum groups. Quantum group isomorphism is not regarded
in this paper.

\section{Quantum complex kinematic groups}

Kinematic groups are  motion groups of the maximal homogeneous four
dimensional (one temporal and three space coordinates)
space--time models \cite{BLL}. All these groups may be obtained from the
real  group $ SO(5;\R) $ by contractions and analytic
continuations \cite{2}.
There are three types of kinematics: nonrelativistic ---
Galilei  $ G(1,3)=SO(5;\iota_{1},\iota_{2},1,1) $
with zero curvature and
Newton  $ N^{\pm}(1,3)=SO(5;j_{1}=1,i;\iota_{2},1,1) $
with positive
and negative curvature, respectively;
relativistic --- Poincare
$ P(1,3)=SO(5;\iota_{1},i,1,1) $
with zero curvature and 
 (anti) de Sitter
$ S^{\pm}(1,3)=SO(5;j_{1}=1,i;i,1,1) $
with (positive)
negative curvature; exotic --- Carroll
$ C^{0}(1,3)=SO(5;\iota_{1},1,1,\iota_{4})$
with zero curvature and $ C^{\pm}(1,3)=SO(5;j_{1}=1,i;1,1,\iota_{4})$
with positive and negative curvature.

The groups $ N^{\pm}(1,3)$ are the real forms of the complex Newton
group $ N(4),$ Poincare group $P(1,3)$ is the real form of the
complex Euclid group $E(4),$ the groups $ C^{\pm}(1,3) $
are the real forms of the complex Carroll group $C(4).$
In this paper the quantum deformations of the complex orthogonal groups
are regarded, therefore whith the help of contractions a quantum analogs
of the complex kinematic groups may be obtained. Possible contractions
of the complex quantum groups $SO_{q}(5;j;\sigma)$ are
described by theorems 1 and 3 for $N=5.$
If  deformation parameter remain unchanged $( J=1),$ then
we have the quantum analogs of  Euclead group $ E_{q}(4), $
  Newton group $ N_{q}(4) $
and   Carroll group $ C_{q}(4). $
If the deformation parameter is transformed under contraction $ z=\iota_{2}v,$
then we have one more quantum deformation of  Newton group $ N_{v}(4), $
which is not isomorphic to the previous one. Two primitive elements of
$ N_{q}(4)$ correspond to the elliptic translation along the temporal axis
$ t $ and to the rotation in the space plane $ \{r_{2},r_{3}\} $ (both
are isomorphic to $ SO(2)$), while primitive elements of $ N_{v}(4)$
correspond to the flat translation along the spatial axis $ r_{2}$
and to  Galilei boost in the space-time plane $ \{t,r_{1}\}$
(both are isomorphic to Galilei group $ SO(2;j_{2}=\iota_{2})=G(1,1)).$
We did not obtain the quantum deformations of the complex Galilei
$G(4)$ and Carroll $C^0(4)$ groups.

According with correspondence principle a new physical theory must include
an old one as a particle case. For space-time theory this principle is
realised as the chain of limit transitions: general relativity pass to
special relativity when space-time curvature tends to zero and special
relativity pass to classical physics when light velosity tends to infinite.
For kinematical groups this correspond to the chain of contractions:
  \begin{equation}
  S^{\pm}(1,3)\stackrel{ K \rightarrow 0}{\longrightarrow}
  P(1,3)\stackrel{c \rightarrow \infty }{\longrightarrow}G(1,3).
    \label{8}
     \end{equation}
As it was mentioned above there is no quantum deformation of the complex
Galilei group, therefore it is not possible to construct the quantum
analog of the full chain of contractions (\ref{8}) even at the level
of complex groups. This means that (at least standart) quantum
deformation of the flat nonrelativistic space-time  do not exist.

\end{document}